\documentclass{article}
\usepackage{amsmath,amssymb,latexsym,amsthm,amscd}

\newcommand\Hh{\mathcal H}

\newcommand\M{\mathcal M}

\newcommand\La{\mathcal L}

\newcommand\bq{\mathbf{Q}}

\newcommand\bk{\mathbf{K}}

\newtheorem{theorem}{Theorem}

\numberwithin{proposition}{section}
\numberwithin{definition}{section}
\numberwithin{corollary}{section} \numberwithin{remark}{section}
\numberwithin{lemma}{section} \numberwithin{equation}{section}

\numberwithin{theorem}{section}

\numberwithin{question}{section} \numberwithin{case}{section}
\numberwithin{example}{section} \numberwithin{conjecture}{section}

\begin{document}

\title{Finiteness of the number of arithmetic groups generated by
reflections in Lobachevsky spaces}
\author{Viacheslav V. Nikulin \footnote{Supported by EPSRC grant
EP/D061997/1}} \maketitle
 
\begin{abstract} After results by the author (1980, 1981) and  
Vinberg (1981),  finiteness of the number of maximal arithmetic groups 
generated by reflections in Lobachevsky spaces was not known in 
dimensions $2\le n\le 9$ only. 

Recently (2005), the finiteness was proved in dimension 2 by 
Long, Maclachlan and Reid, and in dimension 3 by Agol.

Here we use these results in dimensions 2 and 3 to prove finiteness in 
all remaining dimensions $4\le n\le 9$. Methods of the author (1980, 1981) 
are more than sufficient to prove this by a very short and very simple 
consideration. 
\end{abstract}

\section{Introduction} \label{introduction} There are three types of  
simply-connected complete Riemannian manifolds of the constant curvature:
spheres, Euclidean spaces and Lobachevsky (hyperbolic)
spaces. Discrete reflection groups (i. e. generated by
reflections in hyperplanes of these spaces) were defined by H.S.M. 
Coxeter. He classified these groups in spheres and
Euclidean spaces. See \cite{Coxeter}.

There are two types of discrete reflection groups with fundamental
domain of finite volume in Lobachevsky spaces: arithmetic and
general. The subject of this paper is finiteness of the number of
maximal arithmetic reflection groups in Lobachevsky spaces. We
assume that the curvature is $-1$.

First, we remind known results about classification of arithmetic
reflection groups in Lobachevsky spaces.

In \cite{Vin1}, \`E.B. Vinberg (1967) gave the criterion of
arithmeticity for discrete reflection groups in Lobachevsky spaces
in terms of their fundamental chambers.  In particular, he introduced
the notion of the ground field of such group. This is a purely
real algebraic number field of a finite degree over $\bq$. An
arithmetic reflection group is a subgroup of finite index of the
automorphism group of a hyperbolic quadratic form over the ring of
integers of this field. Thus, there are two important integer
parameters related to the arithmetic reflection group: the dimension $n\ge 2$ 
of the Lobachevsky space, and the degree $N\ge 1$ of the ground field over 
$\bq$.

In \cite{Nik1}, the author proved (1980) that the number of maximal
arithmetic reflection groups is finite for the fixed parameters $n$ and $N$.
Moreover, the number of ground fields of the fixed degree $N$ is also finite.

In \cite{Nik2}, the author proved (1981) that there exists a constant $N_0$
such that, the dimension $n\le 9$
if the degree $N\ge N_0$. Thus, 
the number of maximal arithmetic reflection groups is finite in each 
dimension $n\ge 10$ of Lobachevsky space.

In \cite{Vin2} and \cite{Vin3}, \`E.B. Vinberg proved (1981) that 
the dimension $n < 30$. Thus, the number of maximal arithmetic reflection 
groups is finite in all dimensions $n\ge 10$. 

Therefore, in the problem of finiteness of the number of arithmetic reflection 
groups in Lobachevsky spaces, only boundedness of the  
degree $N$ of ground fields in dimensions $2\le n\le 9$ 
of Lobachevsky spaces remained open.  

See also reports \cite{Vin4} by \`E.B. Vinberg (1983) and \cite{Nik3} by the
author (1986) in International Congresses of Mathematicians about these 
results.

During almost 25 years, there were no new general results in this domain.  
But, recently, some new important general results were obtained. 

In \cite{LMR}, D.D. Long, C. Maclachlan, A.W. Reid proved (2005)
that the number of maximal arithmetic reflection groups is finite in
dimension  $n=2$. More generally, they proved finiteness of the 
number of maximal arithmetic
Fuchsian groups of genus $0$.

In \cite{Agol}, I. Agol proved (2005) that the number of maximal
arithmetic reflection groups is finite in dimension $n=3$.

Thus, after all these results, finiteness of the number of maximal
arithmetic reflection groups is not proved in dimensions
$4\le n\le 9$ only.

The aim of this short paper is to show that the authors methods developed 
in  \cite{Nik1} and \cite{Nik2} are more than sufficient 
to deduce this from finiteness results in dimensions $n=2$ and
$n=3$. We don't use at all the very important part of the methods  
which is necessary to bound the dimension of the Lobachevsky space. 

More exactly, 
we show that if the number of maximal arithmetic reflection groups is 
finite in dimensions $2$ and $3$, then from results of \cite{Nik1} and 
\cite{Nik2}, it follows that the degree of ground fields  in dimensions 
$4\le n\le 9$ is bounded. By \cite{Nik1}, it follows finiteness of 
the number of maximal arithmetic reflection groups in these dimensions.

Finally, this completes the proof of finiteness of the number of 
maximal arithmetic reflection groups in Lobachevsky spaces all together. 

\section{Finiteness of the number of maximal arithmetic reflection groups
in Lobachevsky spaces of
dimensions $4\le n\le  9$} \label{sec1}

In the proof below, we work with the Klein model of Lobachevsky space $\La$ 
of dimension $n$ which is related to 
a real hyperbolic form of the signature $(1,n)$. Then a finite 
convex polyhedron $\M$ in $\La$ is defined by the set $P(\M)$ of 
vectors with the square $-2$ of this form which are perpendicular 
to hyperplanes $\Hh_\delta$, $\delta\in P(\M)$, of faces of 
the highest (i. e. of the codimension one) dimension of the polyhedron 
$\M$ and directed outwards of $\M$. Then the number 
$\delta\cdot \delta^\prime$ for $\delta,\delta^\prime\in P(\M)$ 
defines the angle or distance between the corresponding 
hyperplanes $\Hh_\delta$, $\Hh_{\delta^\prime}$. In particular, 
these hyperplanes are perpendicular if and 
only if $\delta\cdot \delta^\prime=0$.   The set $P(\M)$ defines 
the Gram matrix $(\delta\cdot \delta^\prime)$, 
$\delta,\delta^\prime \in P(\M)$, 
of the polyhedron $\M$ and the equivalent Gram diagram $\Gamma(P(\M))$. 
It has vertices $P(\M)$, and two different vertices 
$\delta,\delta^\prime$ are connected by the edge of the weight 
$\delta\cdot \delta^\prime$ if $\delta\cdot \delta^\prime \not=0$. 
This defines connected components of $P(\M)$ and of its subsets. See details 
in \cite{Vin1} and \cite{Nik1}, \cite{Nik2}.  

The aim of this paper is to prove the following theorem. 

\begin{theorem}
\label{maintheorem} In Lobachevsky spaces of dimensions $4\le n\le 9$, 
the degree of ground fields of arithmetic reflection 
groups is bounded.
\end{theorem}

\begin{proof} We consider induction by the dimension $n$ of Lobachevsky spaces.

For $n=2$, the number of maximal arithmetic reflection groups is finite by 
\cite{LMR}. Then, obviously, the degree of ground fields of arithmetic 
reflection groups is bounded in this dimension. 

For $n=3$, the number of maximal arithmetic reflection groups is finite by 
\cite{Agol}, and the degree of ground fields of arithmetic reflection 
groups is also bounded. 

Let us assume that the degree of ground fields of arithmetic reflection groups 
is bounded in all dimensions $\le n-1$ of Lobachevsky
spaces where  $n\ge 4$. Let us prove this for $n$.

Let $G$ be an arithmetic reflection group in Lobachevsky space of
dimension $n$. Let $\M$ be a fundamental chamber of $G$. It is
well-known that $\M$ is bounded (compact) if the ground field is different
from $\bq$. Thus, we can assume that $\M$ is bounded. Let $P(\M)$
be the set of all vectors with square $-2$ which are perpendicular
to codimension one faces of $\M$ and directed outwards. 
We remind that $\M$ has 
acute angles, and then $\delta\cdot \delta^\prime\ge 0$ for different 
$\delta,\delta^\prime \in P(\M)$. 

Following \cite{Nik1} and \cite{Nik2}, let us take a
codimension one face $\Hh_e\cap \M$, $e\in P(\M)$, of $\M$ which
has minimality $14$  (it defines the narrow part of $\M$) (see
Sect. 2 in  \S 4 of \cite{Nik2}). We denote by $P(\M,e)$ the set of
all $\delta \in P(\M)$ such that the hyperplane $\Hh_\delta$
intersects the hyperplane $\Hh_e$ (it is well-known that then 
$\Hh_e\cap \Hh_\delta \cap \M$ is a codimension two face of $\M$
if $\delta\not=e$). The minimality $14$ of the face $\Hh_e\cap \M$ means 
that  we have 
\begin{equation}
\label{narrowcondition}
\delta \cdot \delta^\prime<14,\ \ \forall\, \delta,\delta^\prime\in P(\M,e).
\end{equation}

If $\Hh_\delta\perp \Hh_e$ (equivalently $\delta\cdot e=0$) for all  
$\delta\in P(\M,e)-\{e\}$,
then, obviously, the face $\Hh_e\cap \M$ is a fundamental chamber
for an arithmetic reflection group in $\Hh_e$ of dimension $n-1$
with the same ground field. By induction hypothesis, it has
bounded degree.

Thus, we can assume that there exists $f\in P(\M,e)-\{e\}$ such
that $f\cdot e>0$ (equivalently $f$ and $e$ are connected by an edge 
in the Gram diagram of $\M$).  
Let us consider the codimension two face
$\Hh_e\cap \Hh_f \cap \M$ of $\M$. Let us denote by $P(\M,e,f)$
the set of all $\delta\in P(\M)$ such that the hyperplane
$\Hh_\delta$ intersects the codimension two subspace $\Hh_e\cap
\Hh_f$ (it is well-known that then $\Hh_e\cap \Hh_f\cap \Hh_\delta\cap
\M$ is a codimension $3$ face of $\M$ if $\delta$ is different
from $e$ and $f$). If $\Hh_\delta\perp \Hh_e\cap \Hh_f$ 
(equivalently $\delta\cdot e=\delta\cdot f=0$) for all 
$\delta\in P(\M,e,f)-\{e,f\}$, then, obviously, the codimension
$2$ face $\Hh_e\cap \Hh_f\cap \M$ of $\M$ is a fundamental chamber
of arithmetic reflection group of dimension $n-2\ge 2$ with the same
ground field. By induction hypothesis, the field has  bounded degree.

Thus, we can assume that there exists $g\in P(\M,e,f)-\{e,f\}$
such that $\Hh_g$ is not perpendicular to $\Hh_e\cap \Hh_f$. This
means that either $g\cdot e>0$ or $g\cdot f >0$. Thus, the Gram 
diagram of $e,f,g$ is connected, and their Gram matrix is negative definite. 

Let us take a (one-dimensional)  edge $r$ in the face $\Hh_e\cap
\Hh_f\cap \M$ of $\M$ which terminates in the hyperplane $\Hh_g$.
Thus, one (of two) vertices of $r$ is contained in $\Hh_g$, but
another one is not contained (equivalently, the edge $r$ is not
contained in $\Hh_g$). Existence of such edge $r$ is  obvious.

Let $P(r)$ be the set of all $\delta\in P(\M)$ such that
$\Hh_\delta$ contains at least one of two vertices of the
one-dimensional edge $r$. Since $r\subset \Hh_e\cap \M$, then 
$P(r)\subset P(\M,e)$, and $P(r)$ then defines 
the edge polyhedron of the minimality $14$, since  \eqref{narrowcondition},  
according to the definition 
from (Sect. 2 of \S2 in \cite{Nik2}). 
This edge polyhedron  defines the ground field
$\bk$ of the group $G$. The elements $e,f,g\in P(r)$ have a 
connected Gram diagram, and their hyperplanes contain the same 
vertex of $r$, but the edge $r$ is not contained in $\Hh_g$, 
by our construction.  
This shows that the hyperbolic (i. e. with hyperbolic Gram matrix) 
connected component of the Gram 
diagram $\Gamma(P(r))$ of this edge polyhedron 
is different from each of diagrams $\Gamma_1$, $\Gamma_2$ and $\Gamma_3$ 
of Theorem 2.3.1 from \cite{Nik2}. By the same Theorem 2.3.1 from \cite{Nik2}, 
then $[\bk:\bq]$ is bounded by the constant $N(14)$ of this theorem.

This finishes the proof.
\end{proof}

As we explained in Introduction, this completes the proof of  
finiteness of the number of maximal arithmetic reflection groups in 
Lobachevsky spaces all together.

\begin{theorem}
The number of maximal arithmetic reflection groups in Lobache\-vsky
spaces of dimensions $n\ge 2$ all together is finite.
\end{theorem}

Their number is finite in each dimension $n\ge 4$ (the author).
They don't exist in dimensions $n\ge 30$ (\`E.B. Vinberg).
Their number is finite in the dimension $n=2$ 
(D.D. Long, C. Maclachlan and A.W. Reid).
Their number is finite in the dimension $n=3$ (I. Agol).

V.V. Nikulin \par Deptm. of Pure Mathem. The University of
Liverpool, Liverpool\par L69 3BX, UK; \vskip1pt Steklov
Mathematical Institute,\par ul. Gubkina 8, Moscow 117966, GSP-1,
Russia

vnikulin@liv.ac.uk \ \ vvnikulin@list.ru

\end{document}